\documentclass[10pt]{article}
\usepackage{amscd,amsfonts,amssymb,amsthm,amscd,amsmath,latexsym}

\makeatletter

\usepackage{pslatex}
\begin{titlepage}

\title{\bf \Huge Homotopy quantum field theory and the index gerbe}
\author{Ulrich Bunke\thanks{G\"ottingen, bunke@uni-math.gwdg.de} }
\end{titlepage}
\begin{document} 
\maketitle


\newcommand{\by}{{\bf y}}
\newcommand{\bE}{{\bf E}}
\newcommand{\Face}{{\rm Face}}
\newcommand{\cDelta}{{\bf \Delta}}
\newcommand{\LIM}{{\rm LIM}}
\newcommand{\diag}{{\rm diag}}
\newcommand{\Proof}{{\it Proof.$\:\:\:\:$}}
 \newcommand{\dist}{{\rm dist}}
\newcommand{\kaaa}{{\frak k}}
\newcommand{\paaa}{{\frak p}}
\newcommand{\vp}{{\varphi}}
\newcommand{\taaa}{{\frak t}}
\newcommand{\haaa}{{\frak h}}
\newcommand{\R}{{\Bbb R}}
\newcommand{\Hh}{{\bf H}}
\newcommand{\Rep}{{\rm Rep}}
\newcommand{\Hb}{{\Bbb H}}
\newcommand{\Q}{{\Bbb Q}}
\newcommand{\str}{{\rm str}}
\newcommand{\Ind}{{\rm ind}}
\newcommand{\triv}{{\rm triv}}
\newcommand{\Z}{{\mathbb{Z}}}
\newcommand{\bD}{{\bf D}}
\newcommand{\bF}{{\bf F}}
\newcommand{\tX}{{\tt X}}
\newcommand{\Cliff}{{\rm Cliff}}
\newcommand{\tY}{{\tt Y}}
\newcommand{\tZ}{{\tt Z}}
\newcommand{\tV}{{\tt V}}
\newcommand{\tR}{{\tt R}}
\newcommand{\Fam}{{\rm Fam}}
\newcommand{\Cusp}{{\rm Cusp}}
\newcommand{\bT}{{\bf T}}
\newcommand{\bK}{{\bf K}}
\newcommand{\bo}{{\bf o}}
\newcommand{\K}{{\Bbb K}}
\newcommand{\tH}{{\tt H}}
\newcommand{\bS}{{\bf S}}
\newcommand{\bB}{{\bf B}}
\newcommand{\tW}{{\tt W}}
\newcommand{\tF}{{\tt F}}
\newcommand{\bA}{{\bf A}}
\newcommand{\bL}{{\bf L}}
 \newcommand{\bom}{{\bf \Omega}}
\newcommand{\bundle}{{bundle}}
\newcommand{\ch}{{\bf ch}}
\newcommand{\ve}{{\varepsilon}}
\newcommand{\C}{{\Bbb C}}
\newcommand{\gen}{{\rm gen}}
\newcommand{\cTop}{{{\cal T}op}}
\newcommand{\bP}{{\bf P}}
\newcommand{\Naaa}{{\bf N}}
\newcommand{\image}{{\rm image}}
\newcommand{\gaaa}{{\frak g}}
\newcommand{\zaaa}{{\frak z}}
\newcommand{\saaa}{{\frak s}}
\newcommand{\laaa}{{\frak l}}
\newcommand{\bN}{{\bf N}}
\newcommand{\stimes}{{\times\hspace{-1mm}\bf |}}
\newcommand{\ausg}{{\rm end}}
\newcommand{\bff}{{\bf f}}
\newcommand{\maaa}{{\frak m}}
\newcommand{\aaaa}{{\frak a}}
\newcommand{\naaa}{{\frak n}}
\newcommand{\brr}{{\bf r}}
\newcommand{\res}{{\rm res}}
\newcommand{\Aut}{{\rm Aut}}
\newcommand{\Pol}{{\rm Pol}}
\newcommand{\Tr}{{\rm Tr}}
\newcommand{\cT}{{\cal T}}
\newcommand{\dom}{{\rm dom}}
\newcommand{\Line}{{\rm Line}}
\newcommand{\db}{{\bar{\partial}}}
\newcommand{\Sf}{{\rm  Sf}}
\newcommand{\g}{{\gaaa}}
\newcommand{\cZ}{{\cal Z}}
\newcommand{\cH}{{\cal H}}
\newcommand{\cM}{{\cal M}}
\newcommand{\interi}{{\rm int}}
\newcommand{\singsupp}{{\rm singsupp}}
\newcommand{\cE}{{\cal E}}
\newcommand{\ccR}{{\cal R}}
\newcommand{\hol}{{\rm hol}}
\newcommand{\cV}{{\cal V}}
\newcommand{\cY}{{\cal Y}}
\newcommand{\cW}{{\cal W}}
\newcommand{\del}{{\rm del}}
\newcommand{\bdel}{{\bf del}}
\newcommand{\cI}{{\cal I}}
\newcommand{\cC}{{\cal C}}
\newcommand{\cK}{{\cal K}}
\newcommand{\cA}{{\cal A}}
\newcommand{\cEp}{{{\cal E}^\prime}}
\newcommand{\cU}{{\cal U}}
\newcommand{\Hom}{{\mbox{\rm Hom}}}
\newcommand{\vol}{{\rm vol}}
\newcommand{\cO}{{\cal O}}
\newcommand{\End}{{\mbox{\rm End}}}
\newcommand{\Ext}{{\mbox{\rm Ext}}}
\newcommand{\rk}{{\rm rank}}
\newcommand{\im}{{\mbox{\rm im}}}
\newcommand{\sign}{{\rm sign}}
\newcommand{\spann}{{\mbox{\rm span}}}
\newcommand{\symm}{{\mbox{\rm symm}}}
\newcommand{\cF}{{\cal F}}
\newcommand{\cD}{{\cal D}}
\newcommand{\bC}{{\bf C}}
\newcommand{\bbeta}{{\bf \eta}}
\newcommand{\bOmega}{{\bf \Omega}}
\newcommand{\bbz}{{\bf z}}
\newcommand{\bc}{{\bf c}}
\newcommand{\bb}{{\bf b}}
\newcommand{\bd}{{\bf  d}}
\newcommand{\Ree}{{\rm Re }}
\newcommand{\Res}{{\mbox{\rm Res}}}
\newcommand{\Imm}{{\rm Im}}
\newcommand{\inter}{{\rm int}}
\newcommand{\clo}{{\rm clo}}
\newcommand{\tg}{{\rm tg}}
\newcommand{\ee}{{\rm e}}
\newcommand{\Li}{{\rm Li}}
\newcommand{\cN}{{\cal N}}
 \newcommand{\conv}{{\rm conv}}
\newcommand{\op}{{\mbox{\rm Op}}}
\newcommand{\tr}{{\mbox{\rm tr}}}
\newcommand{\cs}{{c_\sigma}}
\newcommand{\ctg}{{\rm ctg}}
\newcommand{\degg}{{\mbox{\rm deg}}}
\newcommand{\Ad}{{\mbox{\rm Ad}}}
\newcommand{\ad}{{\mbox{\rm ad}}}
\newcommand{\codim}{{\rm codim}}
\newcommand{\Gr}{{\mathrm{Gr}}}
\newcommand{\coker}{{\rm coker}}
\newcommand{\id}{{\mbox{\rm id}}}
\newcommand{\ord}{{\rm ord}}
\newcommand{\nat}{{\Bbb  N}}
\newcommand{\supp}{{\rm supp}}
\newcommand{\sing}{{\mbox{\rm sing}}}
\newcommand{\spec}{{\mbox{\rm spec}}}
\newcommand{\Ann}{{\mbox{\rm Ann}}}
\newcommand{\aca}{{\aaaa_\C^\ast}}
\newcommand{\acag}{{\aaaa_{\C,good}^\ast}}
\newcommand{\acage}{{\aaaa_{\C,good}^{\ast,extended}}}
\newcommand{\tck}{{\tilde{\ck}}}
\newcommand{\tnk}{{\tilde{\ck}_0}}
\newcommand{\ceep}{{{\cal E}(E)^\prime}}
 \newcommand{\ncE}{{{}^\naaa\cE}}
 \newcommand{\Or}{{\rm Or }}
\newcommand{\Diff}{{\cal D}iff}
\newcommand{\cB}{{\cal B}}
\newcommand{\hc}{{{\cal HC}(\gaaa,K)}}
\newcommand{\hcma}{{{\cal HC}(\maaa_P\oplus\aaaa_P,K_P)}}
\def\imath{{\rm i}}
\newcommand{\vsl}{{V_{\sigma_\lambda}}}
\newcommand{\czg}{{\cZ(\gaaa)}}
\newcommand{\csl}{{\chi_{\sigma,\lambda}}}
\newcommand{\cR}{{\cal R}}
\def\hB{\hspace*{\fill}$\Box$ \newline\noindent}
\newcommand{\varho}{\varrho}
\newcommand{\ind}{{\rm index}}
\newcommand{\Indu}{{\rm Ind}}
\newcommand{\Fin}{{\mbox{\rm Fin}}}
\newcommand{\cS}{{S}}
\newcommand{\orig}{{\cal O}}
\def\hB{\hspace*{\fill}$\Box$ \\[0.5cm]\noindent}
\newcommand{\cL}{{\cal L}}
 \newcommand{\cG}{{\cal G}}
\newcommand{\npci}{{\naaa_P\hspace{-1.5mm}-\hspace{-2mm}\mbox{\rm coinv}}}
\newcommand{\pki}{{(\paaa,K_P)\hspace{-1.5mm}-\hspace{-2mm}\mbox{\rm inv}}}
\newcommand{\mki}{{(\maaa_P\oplus \aaaa_P, K_P)\hspace{-1.5mm}-\hspace{-2mm}\mbox{\rm inv}}}
\newcommand{\Mat}{{\rm Mat}}
\newcommand{\npi}{{\naaa_P\hspace{-1.5mm}-\hspace{-2mm}\mbox{\rm inv}}}
\newcommand{\ngp}{{N_\Gamma(\pi)}}
\newcommand{\gbg}{{\Gamma\backslash G}}
\newcommand{\gkm}{{ Mod(\gaaa,K) }}
\newcommand{\ggkm}{{  (\gaaa,K) }}
\newcommand{\pkm}{{ Mod(\paaa,K_P)}}
\newcommand{\ppkm}{{  (\paaa,K_P)}}
\newcommand{\makm}{{Mod(\maaa_P\oplus\aaaa_P,K_P)}}
\newcommand{\mmakm}{{ (\maaa_P\oplus\aaaa_P,K_P)}}
\newcommand{\cP}{{\cal P}}
\newcommand{\gm}{{Mod(G)}}
\newcommand{\gk}{{\Gamma_K}}
\newcommand{\La}{{\cal L}}
\newcommand{\cug}{{\cU(\gaaa)}}
\newcommand{\cuk}{{\cU(\kaaa)}}
\newcommand{\dc}{{C^{-\infty}_c(G) }}
\newcommand{\gdk}{{\gaaa/\kaaa}}
\newcommand{\dgkm}{{ D^+(\gaaa,K)-\mbox{\rm mod}}}
\newcommand{\dgm}{{D^+G-\mbox{\rm mod}}}
\newcommand{\vect}{\mathrm{vect}}
 \newcommand{\cig}{{C^{ \infty}(G)_{K} }}
\newcommand{\gami}{{\Gamma\hspace{-1.5mm}-\hspace{-2mm}\mbox{\rm inv}}}
\newcommand{\cQ}{{\cal Q}}
\newcommand{\mmap}{{Mod(M_PA_P)}}
\newcommand{\bbbz}{{\bf Z}}
 \newcommand{\cX}{{\cal X}}
\newcommand{\bH}{{\bf H}}
\newcommand{\pr}{{\rm pr}}
\newcommand{\bX}{{\bf X}}
\newcommand{\bY}{{\bf Y}}
\newcommand{\bZ}{{\bf Z}}
\newcommand{\bz}{{\bf z}}
\newcommand{\bkappa}{{\bf \kappa}}
\newcommand{\ev}{{\rm ev}}
\newcommand{\bV}{{\bf V}}
\newcommand{\Gerbe}{{\rm Gerbe}}
\newcommand{\gerbe}{{\rm gerbe}}
\newcommand{\hA}{{\bf \hat A}}

\newtheorem{prop}{Proposition}[section]
\newtheorem{lem}[prop]{Lemma}
\newtheorem{ddd}[prop]{Definition}
\newtheorem{theorem}[prop]{Theorem}
\newtheorem{kor}[prop]{Corollary}
\newtheorem{ass}[prop]{Assumption}
\newtheorem{con}[prop]{Conjecture}
\newtheorem{prob}[prop]{Problem}
\newtheorem{fact}[prop]{Fact}

\tableofcontents

\parskip1ex

\section{Introduction}

In the recent paper \cite{turnerwillerton02} Turner and Willerton study the relation
between gerbes with connection and thin-invariant rank-one  field theories on a space $B$. Their main result is that
if $H_1(B,\Z)$ is torsion-free, then there is a one-to-one correspondence of 
gerbes with connection on $B$ and isomorphism classes of thin-invariant rank-one field theories on $B$.
The gerbe corresponding to a thin-invariant rank-one field theory is given by the holonomy of the field theory.

In \cite{lott01} Lott associates to an odd dimensional geometric family on $B$ (see \cite{bunke011} and Section \ref{indexgerbe} for this notion) 
an index gerbe. 

The goal of the present paper is to analyze the construction of Turner and Willerton  \cite{turnerwillerton02} in the case of the index gerbe. Starting from the geometric family we will construct a thin-invariant rank-one field theory on $B$ such that its holonomy is that of the
index gerbe. Our construction works without any assumption on the manifold $B$.
 
If $H_1(B,\Z)$ is torsion-free, then in view of  \cite{turnerwillerton02}, Thm. 3.5, the field theory which we construct in the present paper coincides with the one obtained by applying the 
construction of  Turner and Willerton  \cite{turnerwillerton02} to the index gerbe.

For general $B$ the thin-invariant field theory associated to a geometric family  seems to encode slightly more structure than the index gerbe of that family. The construction of the field theory
given in the present paper induces a construction of the index gerbe which is independent of the
constructions in \cite{lott01} and \cite{bunke02}.

\section{The index gerbe}\label{indexgerbe}

Let $\Gerbe(B)$ denote the group of gerbes with connection on $B$.
We refer to Hitchin \cite{hitchin99} for a nice introduction to the geometric picture of gerbes.
Lott's index gerbe is constructed in \cite{lott01} using  Hitchin's geometric picture. In the present paper we prefer to represent gerbes by Deligne cohomology classes or Cheeger-Simons differential characters.
 
Let  $\cA_B^q$ denote the sheaf of real smooth $q$-forms on $B$, $q\in\nat_0$.
We consider the complex of sheaves on $B$
$$\cK(2,\Z)_B: 0\rightarrow \underline{\Z}_B\rightarrow \cA_B^0\rightarrow \cA_B^1\rightarrow \cA_B^2\rightarrow 0\ ,$$
where the constant sheaf $\underline{\Z}_B$ sits in degree $-1$.
The third Deligne cohomology of $B$ is by definition the
second \v{C}ech hyper cohomology group of this complex
$$H^3_{Del}(B):=\check{\bH}^2(B,\cK(2,\Z)_B)\ .$$
There is a natural isomorphism $\Gerbe(B)\cong H_{Del}^3(B)$ (see 
Brylinski \cite{brylinski93} and \cite{bunke011}).

Holonomy provides a natural isomorphism 
$H:H^3_{Del}(B)\stackrel{\sim}{\rightarrow} \hat H^2(B,U(1))$,
where $\hat H^2(B,U(1))$ is the group of Cheeger-Simons differential characters
with values in $U(1)$ (see \cite{cheegersimons83}). An element $\phi\in \hat H^2(B,U(1))$ 
is a homomorphism from the group of smooth cycles in $B$ to $U(1)$
such that there exists a closed  form $R^\phi\in \cA_B^3(B,d=0)$ with the property that
$$\phi(\partial C)=\exp\left(2\pi\imath \int_C R^\phi\right)$$ for all $3$-chains $C$.
The form $R^\phi$ is called the curvature of $\phi$.
A construction of 
$H:H_{Del}^3(B)\stackrel{\sim}{\rightarrow}  \hat H^2(B,U(1))$ was given in 
\cite{bunke02}, Sec. 6.1.
Thus we can represent gerbes in a third way, namely by  Cheeger-Simons differential characters.

There  is a natural homomorphism $R:H^3_{Del}(B)\rightarrow \cA^3_B(B,d=0)$, $x\mapsto R^x$,
such that $R^x$ is the curvature of the Cheeger-Simons differential character $H(x)$.
We call $R^x$ the curvature of $x$.

An odd dimensional  geometric family $\cE_{geom}$ over $B$ is given by (see \cite{bunke011}, Sec.1.1)
 \begin{itemize}
\item
a fibre bundle $\pi:E\rightarrow B$ with closed odd dimensional fibres, 
\item  a vertical orientation, a vertical spin structure, and a complex vector bundle $V$ over $E$
\item  a vertical Riemannian metric and  a horizontal distribution for $\pi:E\rightarrow B$, a hermitian  metric $h^V$ and a metric  connection $\nabla^V$ on $V$
(the geometric structures). 
\end{itemize} 
By $\gerbe(\cE_{geom})\in \Gerbe(B)$ we denote the index gerbe constructed by Lott \cite{lott01}.
Let $\ind(\cE)\in K^1(B)$ be the index of the family of selfadjoint twisted Dirac operators
defined by $\cE_{geom}$. We omit the subscript ${}_{geom}$ since this index is independent of the choice of the geometric structures. Let $\ch_1(\ind(\cE))\in H^1_{dR}(B)$ be the degree-one component
of the Chern character in the de Rham cohomology of $B$. 
Under the assumption that  $\ch_1(\ind(\cE))=0$  in \cite{bunke02} we constructed
a natural class $\ind_{Del}^3(\cE_{geom})\in H_{Del}^3(B)$ such that
$R^{\ind_{Del}^3(\cE_{geom})}=\Omega^3(\cE_{geom})$, where
$\Omega(\cE_{geom})=\int_{E/B}\hA(\nabla^{T^v\pi})\ch(\nabla^V)$ is the local index form
(see \cite{bunke02}, Sec. 3.2 for definitions).
The class $\ind_{Del}^3(\cE_{geom})$ coincides with $\gerbe(\cE_{geom})$ under the natural isomorphism
$\Gerbe(B)\cong \hat H^2(B,U(1))$.

The construction of the thin-invariant field theory given below will provide an independent construction of the index gerbe as an element in $\hat H^2(B,U(1))$.

\section{Thin invariant field theories}

The thin homotopy category $\cT_B$ of $B$ is defined e.g. in \cite{turnerwillerton02}.
Let $S^1\subset\C$ be the unit circle with a fixed orientation. By $\bar S^1$ we denote the circle with the opposite orientation.
By $S_{m,n}$, $n,m\in\nat_0$, we denote the disjoint union of $m$ copies of $S^1$ and $n$ copies of $\bar S^1$.
An object of $\cT_B$ is a smooth map $\gamma : S_{m,0}\rightarrow B$.

Let $I:=[0,1]$ denote the unit interval with the canonical orientation. We consider the oriented  surface $C_{m,n}:=I\times S_{m,n}$ 
which comes with a natural projection $\pr:C_{m,n}\rightarrow S_{m,n}$.
By a surface $\Sigma$ we mean a compact smooth oriented surface with boundary  $\partial \Sigma$ together with an orientation preserving  collar
$i_{\Sigma}:C_{m,n}\hookrightarrow \Sigma$ such that $i_{\Sigma}(\{0\}\times S_{m,n})=\partial \Sigma$. Here $m,n\in\nat_0$.
The image of $S_{m,0}$ (resp. $S_{0.n}$) is called the ingoing (resp. outgoing ) boundary of $\Sigma$ and will be denoted by $\partial^{in}
\Sigma$ (resp. $\partial^{out}\Sigma$). By $I_1^{in}(\Sigma)$ and $I_1^{out}(\Sigma)$ we denote the set of connected components
of $\partial^{in}\Sigma$ and $\partial^{out}\Sigma$.

A $B$-surface is a smooth map $g:\Sigma\rightarrow B$ such that $g\circ i_{\Sigma}:C_{m,n}\rightarrow B$ factors over the projection
$\pr$. In particular,  a $B$-surface determines an ingoing object $\gamma_{in}:S_{m,0}\rightarrow B$ and an outgoing object $\gamma_{out}:S_{n,0}\rightarrow B$
such that $g\circ i_{\Sigma} = (\gamma\cup \gamma^\prime)\circ \pr :C_{m,n}\rightarrow B$.

By a three manifold $X$ we mean a compact three manifold with corners of codimension at most two and fixed collars.
Let $I_1(X)$ be the set of boundary faces $\partial_i X$, $i\in I_1(X)$.
Then $\partial_i X$, $i\in I_1(X)$,  should be a surface with collared boundary as above.
The collars $i_{i,X}:I\times \partial_i X\hookrightarrow X$ are part of the data.
Let $I_2(X)$ denote the set of faces of codimension two. Let $k\in I_2(X)$ and $\partial_i X$ and $\partial_j X$,
$i,j\in I_1(X)$, be the two boundary faces which meet at $\partial_k X$.
Then we require that $i_{i,X}\circ (\id_I\times i_{\partial_i X})=i_{j,X}\circ (\id_I\times i_{\partial_j X})\circ (\sigma\times\id_{\partial_k X})$ as maps from
$I^2\times \partial_k X\hookrightarrow X$, where $\sigma$ interchanges the components in $I^2$. A $B$-three manifold is a smooth map
$h:X\rightarrow B$ such that the compositions $h\circ i_{i,X}:  I\times \partial_i X \rightarrow B$
factor over the projections $ I\times \partial_i X\rightarrow \partial_iX$ and $B$-surfaces
$g_i:\partial_i X\rightarrow B$ for all $i\in I_1(X)$.

Let $\gamma_{in} :S_{m,0}\rightarrow B$ and $\gamma_{out}:S_{n,0}\rightarrow B$ be two objects of
$\cT_B$. A morphism $\gamma_{in}\rightarrow\gamma_{out}$ is represented by a $B$-surface
with ingoing object $\gamma_{in}$ and outgoing object $\gamma_{out}$.
Two such $B$-surfaces $g$ and $g^\prime$ define the same morphism, if they are thin homotopic,
i.e. there is a homotopy $h:I\times \Sigma\rightarrow B$ from $g$ to $g^\prime$
such that $\rk(dh)\le 2$ everywhere. Here $I\times\Sigma$ is a three manifold in the natural way, and $h:I\times \Sigma\rightarrow B$
is required to be a $B$-three manifold.

Composition of morphisms is induced by glueing.
The category $\cT_B$ has a product induced by disjoint union.

Let $\vect_1$ denote the monoidal category of one dimensional complex vector spaces.
A rank-one thin-invariant field theory (compare \cite{turnerwillerton02}, Def. 3.2) is a monoidal functor
$E:\cT_B\rightarrow\vect_1$ satisfying the following condition: 
There is a closed form $R^E\in\cA^3_B(B,d=0)$ (the curvature of $E$) such that
if $h:C\rightarrow B$ is a $B$-three manifold  with boundary $g:\Sigma\rightarrow B$ such that $\Sigma$ is closed,
then $$E([g])=\exp\left(2\pi\imath \int_{C} h^*R^E\right)\ .$$
Here $[g]$ denotes the class of $g$ in the endomorphisms of the empty
object $\emptyset$ of $\cT_B$, and $E([g])\in \C$ under the natural identification of the linear
endomorphisms of $E(\emptyset)\cong \C$ with $\C$.

To any thin-invariant rank-one field theory $E$ with curvature $R^E$ on $B$ we can associate
a gerbe on $B$ with the same curvature. If $g:\Sigma\rightarrow B$ is a closed $B$-surface, then let $[g,\Sigma]$ denote
the corresponding cycle on $B$. 
The Cheeger-Simons differential character $\phi\in \hat H^2(B,U(1))$ corresponding to the gerbe 
satisfies $\phi([g,\Sigma]):=E([g])$.

\section{Taming}

Let $\cE_{geom}$ be an odd dimensional  geometric family over $B$.
By $D(\cE_{geom}):=(D_b)_{b\in B}$ we denote  the family of twisted Dirac operators given by $\cE_{geom}$.
 
Let $g:M\rightarrow B$ be any smooth map from some manifold $M$ to $B$.
A taming of the family $g^*\cE_{geom}$ is given by a smooth family of selfadjoint smoothing operators
$(Q_m)_{m\in M}$ such that $D_m+Q_m$ is invertible for all $m\in M$.
It  was shown in \cite{melrosepiazza97}, that $g^*\cE_{geom}$ admits a taming exactly if
$\ind(g^*\cE)=0$ holds in $K^1(M)$. 

Let $N\subset M$ and $\gamma:N\rightarrow B$ be the restriction of $g$ to $N$.
Assume that we already have a taming $(\gamma^*\cE_{geom})_t$. 
Then the obstruction against extending the taming  to $M$ is the index element
$\ind(g^*\cE,(\gamma^*\cE_{geom})_t)\in K^1(M,N)$.
In order to construct this element we let $C_N(M)$ be the algebra of
continuous functions on $M$ which vanish on $N$. Let $U$ be some neighborhood of
$N$ to which we can extend the family $(Q_n)_{n\in N}$ giving the taming $(\gamma^*\cE_{geom})_t$. Let $\chi\in C^\infty(M)$ be such that
$\chi_{|N}=1$ and $\chi_{|M\setminus U}=0$. Then the family of operators
$(D_m+\chi(m)Q_m)_{m\in M}$ defines the  element of $KK_1(\C,C_N(M))\cong K^1(M,N)$.

\section{A thin-invariant rank-one field theory associated to an odd dimensional geometric family}

\begin{theorem}
If $\cE_{geom}$ is an odd dimensional geometric family on a connected manifold $B$ such that $\ch_1(\cE_{geom})=0$,
then there is a natural rank-one thin invariant field theory $E$ with curvature $R^E=\Omega^3(\cE_{geom})$
such that the associated gerbe is the index gerbe $gerbe(\cE_{geom})$.
\end{theorem}
\Proof 
We fix a base point $b_0$ and a taming $\cE_{b_0,t}$ of the fibre of $\cE_{geom}$ over $b_0$.

Let $\gamma:S^1\rightarrow B$ be a smooth map.
Since $S^1$ is one dimensional we have an inclusion $\ch_1:K^1(S^1)\hookrightarrow  H^1_{dR}(S^1)$.
Since $\ch_1(\ind(\gamma^*\cE))=\gamma^* \ch_1(\ind(\cE))=0$ the family
$\gamma^*\cE_{geom}$ admits a taming $(\gamma^*\cE_{geom})_t$.

The set of homotopy classes of tamings of $\gamma^*\cE_{geom}$
is a torsor over $K^1(S^1\times I, \partial (S^1\times I))\cong \Z$.
We employ the choice $\cE_{b_0,t}$ in order to distinguish one of these classes.
Let $M:=S^1\cup_{1} I$, where we identify $1\in S^1$ and $1\in I$.
We extend $\gamma$ to $\tilde \gamma:M\rightarrow B$ such that $\tilde \gamma$ maps
$0\in I$ to $b_0$. Note that $K^1(S^1\cup_1  I, S^1\cup \{0\})\cong \Z$.
The distinguished homotopy class of tamings on
$\gamma^*\cE_{geom}$ is now characterized by the property that it can be connected with $\cE_{b_0,t}$
along the path $\tilde\gamma_{|I}$. Because of our assumption $\ch_1(\ind(\cE))=0$
this class is independent of the choice of the path $\tilde\gamma_{|I}$.

We now construct the space $E(\gamma)$.
Let $\gamma:S^1\rightarrow B$ be a smooth map.
We consider the set $\tilde E(\gamma)$ of pairs
$(\lambda, (\gamma^*\cE_{geom})_t)$, where
$\lambda\in\C$ and $(\gamma^*\cE_{geom})_t$ is a taming of  $\gamma^* \cE_{geom}$
in the distinguished homotopy class.

Let $(\gamma^*\cE_{geom})_t$ and $(\gamma^*\cE_{geom})_{t^\prime}$ be two tamings in the distinguished class.
We consider $\tilde \gamma:I\times S^1\stackrel{\pr}{\rightarrow }S^1\stackrel{\gamma}{\rightarrow} B$.
Then we have a taming of $\tilde\gamma^* \cE_{geom}$ over 
$\{0\}\times S^1$ given by  $(\gamma^*\cE_{geom})_t$ and over $\{1\}\times S^1$ given  by  $(\gamma^*\cE_{geom})_{t^\prime}$.
Since these belong to the same homotopy class we can extend this taming to a taming
$(\tilde \gamma ^* \cE_{geom})_t$. Any two of these  extensions are homotopic since $K^1(I\times I \times S^1,\partial (I\times I\times S^1))=0$.

By $\eta^2((\tilde \gamma ^* \cE_{geom})_t)\in \cA^2_{I\times S^1}(I\times S^1)$ we denote the $2$-form component of 
the $\eta$ form as defined in \cite{bunke02}, Sec. 3.3.
We claim that $\int_{I\times S^1} \eta^2((\tilde \gamma ^* \cE_{geom})_t)$ is independent
of the choice of the taming. We consider a homotopy between two choices $(\tilde \gamma ^* \cE_{geom})_t$, $(\tilde \gamma^* \cE_{geom})_{t^\prime}$, i.e.  $\hat \gamma : I\times I\times S^1\stackrel{\pr}{\rightarrow} S^1 \stackrel{\gamma}{\rightarrow} B$ and
a taming $(\hat\gamma^*\cE_{geom})_t$. 
Then we have by \cite{bunke02}, Prop. 3.2,  that $d\eta^2((\hat\gamma ^* \cE_{geom})_t)=\Omega^3(\hat\gamma^*\cE_{geom})$.
Moreover $\Omega^3(\hat\gamma^*\cE_{geom})=0$ since $\hat\gamma$ factors over a one dimensional manifold.
Therefore by Stoke's Lemma 
\begin{eqnarray}
0&=&\int_{\{1\}\times I\times S^1} \eta^2((\hat \gamma ^* \cE_{geom})_t) - \int_{\{0\}\times I\times S^1} \eta^2((\hat \gamma ^* \cE_{geom})_t)\label{t1}\\
&&-\int_{I\times\{1\}\times S^1} \eta^2((\hat \gamma ^* \cE_{geom})_t) + \int_{I\times \{0\} \times S^1} \eta^2((\hat \gamma ^* \cE_{geom})_t)\label{t2}\ .
\end{eqnarray}
Here (\ref{t1}) is equal to 
$$\int_{I\times S^1} \eta^2((\tilde \gamma ^* \cE_{geom})_{t^\prime})- \int_{I\times S^1}\eta^2((\tilde \gamma ^* \cE_{geom})_t) \ .$$
Moreover, $((\hat \gamma ^* \cE_{geom})_t)_{|I\times\{0\}\times S^1} =\pr^* (\gamma^*\cE_{geom})_t$, where
$\pr: I\times\{0\}\times S^1\rightarrow S^1$.
Thus for dimensional reasons  $\eta^2((\hat \gamma ^* \cE_{geom})_t)_{|I\times\{0\}\times S^1}=\pr^*\eta^2((\gamma^*\cE_{geom})_t)=0$.
In a similar manner we have $\eta^2((\hat \gamma ^* \cE_{geom})_t)_{|I\times\{1\}\times S^1}=0$.
Thus the  terms in (\ref{t2}) vanish.
This finishes the proof of the claim.

We now define  an equivalence relation on $\tilde E(\gamma)$ such that
$(\lambda, (\gamma^\cE_{geom})_t)\sim (\lambda^\prime, (\gamma^\cE_{geom})_{t^\prime})$ if
$$\lambda^\prime = \lambda \exp\left(-2\pi\imath \int_{I\times S^1} \eta^2((\tilde \gamma ^* \cE_{geom})_t)\right)\ ,$$
where $(\tilde \gamma ^* \cE_{geom})_t$ is any homotopy between $(\gamma^*\cE_{geom})_t$
and $(\gamma^*\cE_{geom})_{t^\prime}$.

We set $E(\gamma):=\tilde E(\gamma)/\sim$. Then $E(\gamma)$ is a one dimensional complex vector space.
If $\gamma=\cup_{i=1,\dots,m}\gamma_i: S_{m,0}\rightarrow B$, then we set 
$E(\gamma):=\otimes_{i=1}^m E(\gamma_i)$.

We now define $E([g])$ for a morphism $[g]$ in $\cT_B$ represented
by a $B$-surface $g:\Sigma\rightarrow B$.

First assume that $\Sigma$ is closed. We again have an inclusion
$\ch_1:K^1(\Sigma)\hookrightarrow H^1_{dR}(\Sigma)$.
Thus $\ind(g^*\cE)=0$ because of our assumption $\ch_1(\ind(\cE))=0$,
and there exists a taming $(g^*\cE_{geom})_t$.
The set of homotopy classes of tamings is parameterized
by $K^1(I\times\Sigma,\partial I\times \Sigma)\cong K^0(\Sigma)\cong \Z\oplus \Z$.
We restrict the choices of tamings $(g^*\cE_{geom})_t$
by the following condition. Let $s\in\Sigma$ be any point.
We consider $M:=\Sigma\cup_{s\sim 1}I$. We choose some extension
$\tilde g: M\rightarrow B$ of $g$ by choosing a path from $b_0$ to $g(s)$.
We will only consider tamings $(g^*\cE_{geom})_t$
which can be connected with $\tilde g_{|\{0\}}^*\cE_{b_0,t}$ along $M$.
This condition is independent of the choice of $s$ and the path.
The remaining choices correspond to the reduced $K$-theory $\tilde K^0(\Sigma)\cong \Z$.

We define 
$$E(g):=\exp\left(2\pi\imath \int_{\Sigma} \eta^2((g^*\cE_{geom})_t)\right) \in \C\cong \End(E(\emptyset))\ .$$

If $\Sigma$ has a boundary, then for
each boundary component $i\in I_1(\Sigma)$ we obtain an object
$\gamma_i : S^1\rightarrow B$. 
We choose tamings $(\gamma_i^*\cE_{geom})_t$, $i\in I_1(\Sigma)$, in the distinguished
components. The obstruction against extending these tamings to $g^*\cE_{geom}$ belongs to
$K^1(\Sigma,\partial \Sigma)$. Since this group is trivial we can 
extend the given taming over the boundaries to a taming $(g^*\cE_{geom})_t$. 
The set of  homotopy classes of these extensions is parameterized by
$K^1(I\times\Sigma,\partial (I\times\Sigma))\cong K^0(\Sigma,\partial \Sigma)\cong \Z$.

Let
$z\in E(\gamma_{in})$ be represented by $$\otimes_{i\in I_1^{in}(\Sigma) }(\lambda_i , (\gamma_i^*\cE_{geom})_t)\ .$$
Then we define 
$$E(g)z:= \exp\left(2\pi\imath \int_{\Sigma} \eta^2((g^*\cE_{geom})_t)\right)
z^\prime\in E(\gamma_{out})\ ,$$
where $z^\prime\in E(\gamma_{out})$ is  represented by $$\otimes_{i\in I_1^{out}(\Sigma)} (\lambda_i , (\gamma_i^*\cE_{geom})_t)\ .$$

 We must show that $E(g)$ only depends on the class $[g]\in \cT_B$.

Let $g$ and $g^\prime$ be thin homotopic by a homotopy $h:I\times \Sigma\rightarrow B$.
Furthermore, let $(h^*\cE_{geom})_t$ be a taming which restricts to the tamings
$(g^*\cE_{geom})_t$ on $\{0\}\times\Sigma$ and $(g^{\prime*}\cE_{geom})_{t^\prime}$ on $\{1\}\times\Sigma$.
Then we have by   \cite{bunke02}, Prop. 3.2., that $d\eta^2((h^*\cE_{geom})_t)=\Omega^3(h^*\cE_{geom})$.
Since the homotopy $h$ is thin we have $ \Omega^3(h^*\cE_{geom})=h^* \Omega^3(\cE_{geom})=0$.
  It follows from Stoke's Lemma that  
$$\int_{\Sigma}
\eta^2((g^*\cE_{geom})_t)-\int_{I\times \partial_{in} \Sigma}\eta^2((h^*\cE_{geom})_t) +
\int_{I\times \partial_{out} \Sigma}\eta^2((h^*\cE_{geom})_t)
=\int_{\Sigma}
\eta^2((g^{\prime*}\cE_{geom})_{t^\prime})\ .$$
This implies that $E(g)=E(g^\prime)$.

Assume that we are given two tamings $(g^*\cE_{geom})_t$ and  $(g^*\cE_{geom})_{t^\prime}$
which coincide over the boundary of $\Sigma$, and which are  not 
homotopic. Let $\tilde \Sigma =\Sigma\cup_{\partial \Sigma}(\bar\Sigma)$ be the surface
obtained by doubling $\Sigma$ along the boundary, and let $\tilde g:\tilde \Sigma\rightarrow B$ be the map induced by
$g$. We consider the  family of cylinders
$I\times \tilde g^*E\rightarrow \tilde \Sigma$ which is boundary tamed by $(g^*\cE_{geom})_t$ on 
$\{0\}\times \tilde g^*E$, and by $(g^*\cE_{geom})_{t^\prime}$ on $\{1\} \times\tilde g^*E$ 
over the copy $\Sigma$, and similarly by $(g^*\cE_{geom})_t$ and again  $(g^*\cE_{geom})_t$ over the other copy $\bar\Sigma$.
This works since the two tamings $(g^*\cE_{geom})_t$ and  $(g^*\cE_{geom})_{t^\prime}$ coincide near $\partial \Sigma$.
The index $\ind((I\times  \tilde g^*\cE_{geom})_{bt})\in K^0(\tilde \Sigma)$ is the obstruction against
extending the boundary taming to a taming (we refer to \cite{bunke02}, Sec. 2.4 and 2.5 for definitions).
It follows from our restrictions on the choice of tamings 
$(g^*\cE_{geom})_t$ and  $(g^*\cE_{geom})_{t^\prime}$ that $\dim \ind((I\times  \tilde g^*\cE_{geom})_{bt})=0$.

Let $\cF_{geom}$ be any geometric family over $\tilde\Sigma$ with closed even dimensional fibres such that
$$\ind(\cF) = -\ind((I\times \tilde g^*\cE_{geom})_{bt})\ .$$ In fact  we can realize $\cF_{geom}$ with zero dimensional fibres
(see \cite{bunke02}, Sec. 3.1). Then $\ind((I\times   \tilde g^*\cE_{geom} + \cF_{geom})_{bt})=0$, and  
 the boundary taming $(I\times  \tilde g^*\cE_{geom} + \cF_{geom})_{bt}$ can be extended to a taming.
Using \cite{bunke02}, Prop. 3.2, (at $*$), we obtain    
 \begin{eqnarray*}
\int_\Sigma\eta^2((g^*\cE_{geom})_{t^\prime})-\int_\Sigma\eta^2((g^*\cE_{geom})_t)&=&
\int_{\tilde \Sigma }\eta^2 (\partial (I\times \tilde g^*\cE_{geom} + \cF_{geom})_{t}) \\
&\stackrel{*}{=}&- \int_{\tilde \Sigma }d \eta^1((I\times \tilde g^*\cE_{geom} + \cF_{geom})_t) +\int_{\tilde\Sigma}\Omega^2(\cF_{geom})\\
&=&<\ch_2(\ind(\cF) ) ,[\tilde \Sigma]>\\
&\in&\Z\ .
\end{eqnarray*}
We conclude that
$$\exp\left(2\pi\imath \int_{\Sigma}
\eta^2((g^*\cE_{geom})_t)\right)= \exp\left(2\pi\imath \int_{\Sigma}
\eta^2((g^{\prime,*}\cE_{geom})^\prime_t)\right)\ .$$
We now have seen that $E(g)$ only depends on the class $[g]\in\cT_B$.

We have constructed for each choice of a taming $\cE_{b_0,t}$
a thin-invariant rank-one field theory. We must now show that
this theory is independent of this choice.

Let $\cE_{b_0,t^\prime}$ be another choice, and let $E^\prime$ be the corresponding
field theory. We must define an equivalence $Q:E\rightarrow E^\prime$.
Let $\gamma:S^1\rightarrow B$ be an object of $\cT_B$.
Let $(\lambda, (\gamma^*\cE_{geom})_t)$ represent
$x\in E(\gamma)$, and let $(\gamma^*\cE_{geom})_{t^\prime}$ be a taming
in the distinguished class for the choice $\cE_{b_0,t^\prime}$.
Then we consider the geometric family $I\times\gamma^*\cE_{geom}$
over $S^1$ with a boundary taming by $ (\gamma^*\cE_{geom})_t$ and $(\gamma^*\cE_{geom})_{t^\prime}$
at $\{0\}\times \cE_{geom}$ and $\{1\}\times \cE_{geom}$.
We consider $\ind((I\times\gamma^*\cE_{geom})_{bt})\in K^0(S^1)\cong \Z$.
Let $\cF$ be a geometric family over $S^1$ with underlying bundle $S^1\rightarrow S^1$ and vector bundle $\C^n\times S^1\rightarrow S^1$,
where we adjust $n$ and the fibre wise orientation such that $\ind(\cF)=-\ind((I\times\gamma^*\cE_{geom})_{bt})$.
Since $\ind((I\times\gamma^*\cE_{geom} \cup\cF_{geom})_{bt})=0$  this boundary taming can be extended
to a taming $(I\times\gamma^*\cE_{geom} \cup\cF_{geom})_t$.

We define 
$$Q(\gamma)(x)=\left(\lambda \exp\left(-2\pi\imath \int_{S^1} \eta^1((I\times\gamma^*\cE_{geom} \cup\cF_{geom})_t)\right) , (\gamma^*\cE_{geom})_{t^\prime}\right)\ .$$
We claim that $Q$ is well-defined.
Let $(\gamma^*\cE_{geom})_{\tilde t}$ and $(\gamma^*\cE_{geom})_{\tilde t^\prime}$ other choices
 of the tamings in the distinguished classes. Let $(\tilde \gamma^*\cE_{geom})_t$ and $(\tilde \gamma^*\cE_{geom})_{t^\prime}$
be corresponding homotopies of tamings from
$(\gamma^*\cE_{geom})_{t}$ to $(\gamma^*\cE_{geom})_{\tilde t}$  and
from $(\gamma^*\cE_{geom})_{t^\prime}$ to $(\gamma^*\cE_{geom})_{\tilde t^\prime}$,
where $\tilde\gamma:I\times S^1\stackrel{\pr}{\rightarrow} S^1 \stackrel{\gamma}{\rightarrow} B$.

Then we obtain a boundary tamed family $(I\times \tilde\gamma^*\cE_{geom})_{bt}$ over $I\times S^1$.
Let $\tilde\cF_{geom}:=\pr^*\cF_{geom}$.
Then we can extend the boundary taming of $I\times\tilde \gamma^*\cE_{geom} \cup\tilde\cF_{geom}$
to a taming $(I\times\tilde \gamma^*\cE_{geom} \cup\tilde\cF_{geom})_{t}$.
By \cite{bunke02}, Prop. 3.2, we have
$$d\eta^1((I\times\tilde \gamma^*\cE_{geom} \cup\tilde\cF_{geom})_{t})=-\eta^2((\tilde \gamma^*\cE_{geom})_{t^\prime})
 +\eta^2((\tilde \gamma^*\cE_{geom})_{t})\ .$$
By Stoke's Lemma
$$-\int_{S^1} \eta^1((I\times \gamma^*\cE_{geom} \cup\cF_{geom})_{\tilde t}) 
+ \int_{S^1} \eta^1((I\times  \gamma^*\cE_{geom} \cup\cF_{geom})_{t})=-\int_{I\times  S^1}\eta^2((\tilde \gamma^*\cE_{geom})_{t^\prime})
+ \int_{I\times  S^1} \eta^2((\tilde \gamma^*\cE_{geom})_{t})\ .$$
An inspection of the definition of $E(\gamma)$ shows that this relation implies that
$Q(\gamma)$ is well-defined independent of the choice of representatives of $E(\gamma)$, $E^\prime(\gamma)$, and of the additional tamings.
We use the monoidal structures of $\cT_B$ and $\vect_1$ in order to extend $Q$ to general objects of $\cT_B$.

It remains to show that $Q$ is natural. Let $g:\Sigma\rightarrow B$ be an $X$-surface
representing a morphism $[g]$ from $\gamma_{in}$ to $\gamma_{out}$.
Let $(g^*\cE_{geom})_t$ and $(g^*\cE_{geom})_{t^\prime}$ be tamings in the distinguished classes with respect  to 
$\cE_{b_0,t}$ and $\cE_{b_0,t^\prime}$.
Then we consider the geometric family $I\times g^*\cE_{geom}$ over $\Sigma$  which is boundary tamed
by $(g^*\cE_{geom})_t$ and  $(g^*\cE_{geom})_{t^\prime}$ at $\{0\}\times g^*\cE_{geom}$ and $\{1\}\times g^*\cE_{geom}$.
We now want to kill $\ind((I\times g^*\cE_{geom})_{bt})\in K^0(\Sigma) $.

Let  $L_0:=\C^{|\dim(\ind((I\times g^*\cE_{geom})_{bt}))|}\times \Sigma \rightarrow \Sigma$ be  a trivial
vector bundle. If $\Sigma$ is closed, then furthermore  let $L_1\rightarrow  \Sigma$ be a line bundle with $c_1(L_1)=-c_1(\ind((I\times g^*\cE_{geom})_{bt}))$. 
We equip $L_1$ with a hermitian metric and a metric connection.
We let $\cF_{0,geom}$ be the geometric family with underlying fibre bundle
$ \Sigma\rightarrow  \Sigma$ and vector bundle $L_0$, where we flip the orientation of the fibres
if $\dim(\ind((I\times g^*\cE_{geom})_{bt}))>0$. If $\Sigma$ is closed, then we let $\cF_{1,geom}$ be the union of geometric families
$\cF_{1a,geom}$ and $\cF_{1b,geom}$.
Here $\cF_{1a,geom}$ has the underlying fibre bundle $\Sigma\rightarrow \Sigma$ and the vector bundle is $L_1$,
and $\cF_{1b,geom}$ has underlying fibre bundle $\Sigma\rightarrow \Sigma$ with flipped orientation and the bundle
is $ \C\times \Sigma\rightarrow \Sigma$. Finally we set $\cF_{geom}:=\cF_{0,geom}$ if $\Sigma$ has a non-trivial boundary, and 
$\cF_{geom}:=\cF_{0,geom}\cup \cF_{1,geom}$ if $\Sigma$ is closed.
Then $\ind((I\times g^*\cE_{geom}\cup \cF_{geom})_{bt})=0$ and 
we can extend the boundary taming  to a taming $(I\times g^*\cE_{geom}\cup \cF_{geom})_t$.
Note that $\Omega^2(I\times g^*\cE_{geom}\cup \cF_{geom})=\Omega^2(\cF_{geom})$
and $\int_{\Sigma} \Omega^2(\cF_{geom})=<c_1(\ind(\cF)),[\Sigma]> \in \Z$, if $\Sigma$ is closed, and
$\int_{\Sigma} \Omega^2(\cF_{geom})=0$, if $\Sigma$ is not closed.
By \cite{bunke02}, Prop. 3.2, we have
$$d\eta^1((I\times g^*\cE_{geom}\cup \cF_{geom})_t) = 
-\eta^2((g^*\cE_{geom})_{t^\prime}) +\eta^2((g^*\cE_{geom})_t) + \Omega^2(\cF_{geom})\ .$$
By Stoke's Lemma
$$\int_{\partial \Sigma}  \eta^1((I\times g^*\cE_{geom}\cup \cF_{geom})_t) = -
\int_{\Sigma} \eta^2((g^*\cE_{geom})_{t^\prime}) +\int_{\Sigma } \eta^2((g^*\cE_{geom})_t)  +\Z\ .$$
This relation implies
$$Q(\gamma) E(g) = E^\prime(g) Q(\gamma)\ .$$

Next we show that the curvature of $E$ is $\Omega^3(\cE_{geom})$.
Let $h:C\rightarrow B$ be a $B$-three manifold with boundary $g:\Sigma\rightarrow B$ such that
$\Sigma$ is closed. If $\partial C=\emptyset$, then
$$\int_{C} \Omega^3(\cE_{geom})=<\ch_3(\ind(\cE)),[C]>\in\Z$$
(see Lott \cite{lott01}, Prop. 8)
so that
$$\exp\left(2\pi\imath \int_{C} \Omega^3(\cE_{geom}) \right)=1=E([g])\ .$$
If $\partial C\not=0$, then $C$ is homotopy equivalent  to a two dimensional space.
In this case vanishing of $\ch_1(\ind(h^*\cE))=h^*\ch_1(\ind(\cE))$ implies that
there exists a taming $(h^* \cE_{geom})_t$. 
We have by \cite{bunke02}, Prop. 3.2,
\begin{eqnarray*}
E([g]) &=&\exp\left(2\pi\imath \int_{\partial C}\eta^2((h^*\cE_{geom})_t)\right)\\
&=&\exp\left(2\pi\imath \int_{C}d\eta^2((h^*\cE_{geom})_t)\right)\\
&=&\exp\left(2\pi\imath \int_{C}\Omega^3(\cE_{geom})\right) \ .
\end{eqnarray*}

Finally we show that the gerbe associated to  the rank-one thin-invariant field theory $E$ is the  index gerbe $\gerbe(\cE_{geom})$.
Let $g:\Sigma\rightarrow B$ be a closed $B$-surface defining the cycle $[g,\Sigma]$.
Let $U$ be a tubular  neighborhood of $g(\Sigma)$.
We choose a taming $\cE_{|U,t}$ which can be connected with $\cE_{b_0,t}$
along some path from $U$ to $b_0$. Then we obtain 
a taming  $(g^*\cE_{geom})_t$. 
Then we have
\begin{eqnarray*}
H(\gerbe(\cE_{geom}))([g,\Sigma])&=&\exp\left(2\pi\imath \int_{\Sigma} g^* \eta^2(\cE_{|U,t})\right)\\& =&
\exp\left(2\pi\imath \int_{\Sigma} \eta^2((g^*\cE_{geom})_t)\right)\\
&=& E([g])\ .
\end{eqnarray*}
\hB
 
\underline{Remark:}
Let $\cE_{geom}$ be an odd-dimensional geometric family over $B$. Let $f:B\rightarrow S^1$
by any classifying map of $-c_1(\ind(\cE))\in H^1(B,\Z)$. Furthermore, let
$\cF_{geom}$ be some  odd dimensional geometric family over $S^1$ such that $c_1(\ind(\cF))=1\in H^1(S^1,\Z)\cong \Z$.
Then we have $\ch_1(\ind(\cE+f^*\cF))=0$. One can show that
the rank-one  thin-invariant field theory associated to $\cE_{geom}+f^*\cF_{geom}$ is independent of the choice
of $f$ and $\cF_{geom}$.

\bibliographystyle{plain}


\begin{thebibliography}{1}

\bibitem{brylinski93}
J.~L. Brylinski.
\newblock {\em Loop Spaces, Characteristic Classes, and Geometric
  Quantization}.
\newblock Birk\"auser, Progress in Math. 107, 1993.

\bibitem{bunke011}
U.~Bunke.
\newblock Transgression of the index gerbe.
\newblock arXiv:math.DG/0109052, 2001.

\bibitem{bunke02}
U.~Bunke.
\newblock Index theory, eta forms, and Deligne cohomology.
\newblock arXiv:math.DG/0201112, 2002.

\bibitem{cheegersimons83}
J.~Cheeger and J.~Simons.
\newblock Differential characters and geometric invariants.
\newblock In {\em LNM1167}, pages 50--80. Springer Verlag, 1985.

\bibitem{hitchin99}
N.~Hitchin.
\newblock Lectures on special lagrangian submanifolds.
\newblock arXiv:math.DG/9907034, 1999.

\bibitem{lott01}
J.~Lott.
\newblock Higher-degree analogs of the determinant line bundle.
\newblock arXiv:math.DG/0106177, 2001.

\bibitem{melrosepiazza97}
R.~B. Melrose and P.~Piazza.
\newblock Families of Dirac operators, boundaries, and the $b$-calculus.
\newblock {\em J. Differential. Geom.}, 46 (1997), 99--180.

\bibitem{turnerwillerton02}
P.~Turner and S.~Willerton.
\newblock Gerbes and homotopy quantum field theory.
\newblock arXiv:math.AT/0201116, 2002.

\end{thebibliography}

\end{document}